\documentclass[a4paper,10pt]{article}
\usepackage{amssymb,amsmath,epsfig,theorem}
\usepackage[english]{babel}

\input xy
\xyoption{all}

\newtheorem{prop}{Proposition}[section]
\newtheorem{theor}[prop]{Theorem}
\newtheorem{lemm}[prop]{Lemma}
\newtheorem{defin}[prop]{Definition}
\newtheorem{corollary}[prop]{Corollary}

\begin{document}

\title{\textbf{On higher dimensional
  Hirzebruch-Jung singularities}} 
\author{Patrick Popescu-Pampu
\\  {\small ppopescu@math.jussieu.fr}}

\date{}
\maketitle

\thispagestyle{empty}
\begin{abstract}{A germ of normal complex analytical surface is called
    a 
    Hirzebruch-Jung singularity if it is analytically isomorphic to the germ
    at the 0-dimensional orbit of an affine toric surface. Two such germs are
    known to be isomorphic if and only if the toric surfaces
    corresponding to them are equivariantly isomorphic. We extend this
    result to higher-dimensional Hirzebruch-Jung singularities, which
    we define to be the germs analytically 
    isomorphic to the germ at the 0-dimensional orbit of an affine toric
    variety determined by a lattice  and a simplicial
    cone of maximal dimension. We deduce a normalization algorithm for 
    quasi-ordinary hypersurface singularities.} 
\end{abstract}


\par\medskip\centerline{\rule{2cm}{0.2mm}}\medskip
\setcounter{section}{0}

\section{Introduction}

In this paper we generalize to arbitrary dimensions the notion of
\textit{Hirzebruch-Jung singularities} and we show how to classify 
them up to analytical isomorphism by combinatorial data. Then we give normal 
forms to these data and we compute these normal forms when the germ is
the normalization of an irreducible quasi-ordinary hypersurface
singularity. 
\medskip

A germ of reduced equidimensional complex analytical space is called
\textit{quasi-ordinary} if there exists a finite morphism from it to a
smooth space of the same dimension, such that the discriminant locus
of the morphism 
is contained in a divisor with normal crossings (see section \ref{gener}).

If the term "quasi-ordinary" seems to appear first in the '60s, in works
of Zariski and Lipman, the study of quasi-ordinary germs goes back at
least to the work \cite{J 08} of Jung on the problem of local
uniformisation of surfaces. For details on it see the first chapter of
\cite{Z 35}. The idea of Jung was to study an arbitrary germ of
surface embedded in $\mathbf{C}^{3}$ by considering a finite linear
projection and an embedded desingularization of the discriminant
curve. By changing the base of the initial projection using this
desingularization morphism, he obtained a surface which is
quasi-ordinary in the neighborhood of any of its points. 

This method was used by Walker \cite{W 35} in order to prove the
existence of a resolution of the singularities of a complex algebraic
surface. This work is considered by Zariski \cite{Z 35} to be the
first rigorous proof of this fact. Hirzebruch \cite{H 53} uses again
Jung's method in order to prove the existence of a desingularization
for complex analytical surfaces which are locally embeddable in
$\mathbf{C}^{3}$. This last restriction was eliminated by Laufer
\cite{L 71}.

An important step in Hirzebruch's method was to consider the normalizations
of the quasi-ordinary germs he arrived at by Jung's method. He gave
explicit constructions of their minimal resolutions by patching affine
planes. Later on, those germs were called "Hirzebruch-Jung
singularities". After Artin's work on rational surface singularities
in the '60s, they were 
seen to be precisely the rational surface singularities which have as
dual resolution graph a segment. This is the definition used in
\cite{BPV 84}. Hirzebruch-Jung singularities are usually classified up
to analytical isomorphism by
an odered pair $(n,q)\in \mathbf{N}^{*}\times \mathbf{N}$ of coprime
numbers with $q <n$. In order to get this classification, Hirzebruch
studied the exceptional divisor of the minimal resolution morphism of
the singularity and introduced the numbers $n,q$ starting from the
self-intersection numbers of its components (see \cite{BPV 84} and
section \ref{class}). It is
also known that this 
classification is topological. For historical details, see
Brieskorn \cite{B 00}.

After the introduction of toric geometry in the '70s, Hirzebruch-Jung
surface singularities were seen to be precisely the germs analytically
isomorphic to the germs of toric surfaces taken at 0-dimensional
orbits (see \cite{O 88} and \cite{F 93}). It is this view-point which
we generalize here. 

 If $\mathcal{W}$ is a lattice and $\sigma$ is a strictly
convex finite rational polyhedral cone in $\mathcal{W}_{\mathbf{R}}:=
\mathcal{W} \otimes \mathbf{R}$, we denote by $\mathcal{M}$ the dual lattice of
$\mathcal{W}$ and by $\check{\sigma}\subset \mathcal{M}_{\mathbf{R}}$
the dual cone of $\sigma$. We denote by 
$\mathcal{Z}(\mathcal{W},
\sigma):=\mathrm{Spec}\:\mathbf{C}[\check{\sigma}\cap \mathcal{M}]$
the affine normal toric variety 
 determined by the pair $(\mathcal{W},\sigma)$.  When $\sigma$ and
$\mathcal{W}_{\mathbf{R}}$ have the same dimension $d$, we say that
$(\mathcal{W}, \sigma)$ is a \textit{maximal pair of
  dimension} $d$.
When $\sigma$ is a simplicial cone, we say that
$(\mathcal{W}, \sigma)$ is a \textit{simplicial pair}. We say that the
simplicial cone $\sigma$ is 
\textit{regular with respect to} $\mathcal{W}$ if it is generated by a
subset of a basis of $\mathcal{W}$. In this case we say also that
$(\mathcal{W}, \sigma)$ is a \textit{regular pair}. 
Two pairs $(\mathcal{W}_{1}, \sigma_{1})$ and
$(\mathcal{W}_{2}, \sigma_{2})$ are called \textit{isomorphic} if
there exists an 
isomorphism of lattices $\phi:\mathcal{W}_{1} \rightarrow
\mathcal{W}_{2}$ sending $\sigma_{1}$ onto $\sigma_{2}$. 

By analogy with the bidimensional case, one can define:
\medskip

\textit{ A germ of irreducible normal complex analytical space of
  arbitrary dimension is called
  \textbf{a Hirzebruch-Jung singularity} if it is analytically isomorphic with 
 the normalization of an $n$-dimensional irreducible quasi-ordinary
germ.} 
\medskip

In \cite{PP 01} (see also \cite{PP 03'} and section \ref{genHJ}) we
showed that such a  normalization is in 
fact analytically isomorphic to the germ at the 0-dimensional orbit of
an affine toric variety defined by a maximal simplicial pair. 
Conversely (see proposition
\ref{equidef}), the germ at the 0-dimensional orbit of a toric variety
defined by a maximal simplicial pair is
quasi-ordinary. This shows that, alternatively, one can define
Hirzebruch-Jung singularities by combinatorial data (see section \ref{genHJ}):
\medskip

\textit{A germ of irreducible normal complex analytical space of
  arbitrary dimension is called
  \textbf{a Hirzebruch-Jung singularity} if it is analytically
  isomorphic with the 
  germ at the 0-dimensio\-nal orbit of an 
affine toric variety defined by a maximal simplicial pair.}
\medskip

It is clear that isomorphic maximal simplicial pairs give rise to
analytically isomorphic  Hirzebruch-Jung singularities. 
Our main theorem (see theorem \ref{mthm}) shows the converse statement:
\medskip

\textit{The analytical type of a Hirzebruch-Jung singularity
  $(\mathcal{Z},0)\simeq(\mathcal{Z}(\mathcal{W}, \sigma),0)$ determines the
  pair $(\mathcal{W}, \sigma)$ up to isomorphism.}
\medskip

In order to prove this result, we make the Riemann extension of the
universal covering map of the smooth part 
of $(\mathcal{Z},0)$ over all
of $(\mathcal{Z},0)$. We call this map $\mu:(\widetilde{\mathcal{Z}},0)
\rightarrow (\mathcal{Z},0)$ \textit{the orbifold map} of
$\mathcal{Z}$ (see section \ref{mainres}). Then we look at the action
$\rho(\mathcal{Z})$ 
of the local fundamental group of $ (\mathcal{Z},0)$ on the Zariski
cotangent space of $(\widetilde{\mathcal{Z}},0)$ and we construct from it
a pair $(W(\rho(\mathcal{Z})), \sigma_{0})$ determined by the
analytical 
type of $\mathcal{Z}$. Theorem \ref{mthm} says that the pairs
$(\mathcal{W}, \sigma)$ and $(W(\rho(\mathcal{Z})), \sigma_{0})$ are
isomorphic. 

We used for the first time orbifold maps in \cite{PP 03'} in order to
get analytical invariants of quasi-ordinary singularities. When we
began to study the problems solved in \cite{PP 03'} and in the present
paper, we tried to use some desingularization morphism of
$(\mathcal{Z},0)$. We could not manage their high non-canonicity, and
so the idea to use instead the orbifold map came as a 
relief. 

We see that, in order to classify up to analytical isomorphism
$n$-dimensional 
Hirzebruch-Jung singularities, one needs only to classify up to
isomorphism the pairs $(\mathcal{W}, \sigma)$, which is a
combinatorial problem. We give normal
forms for such pairs once an ordering of the edges of $\sigma$ is
chosen (proposition \ref{normcone}). We define \textit{the type} of a
Hirzebruch-Jung singularity $(\mathcal{Z},0)$ to be one of the normal
forms associated to the pair $(W(\rho(\mathcal{Z})), \sigma_{0})$
(definition \ref{type}). 

In section \ref{sectalg} we give an \textit{algorithm of
  normalization} of an irreducible quasi-ordinary
hypersurface singularity (proposition \ref{algfond}). More precisely, we
compute the type of the normalization, the ordering being the one
determined by the choice of the ambient coordinates of the starting
quasi-ordinary singularity. The algorithm starts from the 
characteristic exponents and constitute a generalization of the
normalization algorithm for
  surfaces that we published in \cite{PP 01} and \cite{PP
  03} (see proposition \ref{2alg}). Incidentally, if
  $(\mathcal{W},\sigma)$ is a 
  maximal regular pair, we compute the normal forms for the
  pairs $(\mathcal{W}', \sigma)$, where $\mathcal{W}'$ is a sublattice
  of finite index of $\mathcal{W}$ defined by a congruence (lemma
  \ref{cal3}). 

Section \ref{example} contains  a tridimensional example of
application of the algorithm. In 
section \ref{class} we restrict our attention to the bidimensional
case and we compare our definition of the type with the classical one. We
conclude by stating in section \ref{quest} some questions about the
topological types of 
Hirzebruch-Jung singularities and about the analytical types of the
germs at the 0-dimensional orbits of general affine toric varieties.

\medskip

\textbf{Acknowledgements: } I am grateful to Cl{\'e}ment Caubel, Mart{\'\i}n
Sombra and 
Bernard Teissier for their useful remarks on a previous version of
this paper.

\section{Generalities on quasi-ordinary germs} \label{gener}

For any point $P$ on a complex analytical space $\mathcal{V}$, we denote by 
$\mathcal{O}_{\mathcal{V},P}$ the local algebra of $\mathcal{V}$ at
$P$.  In the sequel we will
denote with the same letter a germ and a sufficiently small
representative of it. It will be deduced from the context if one deals
with one or the other notion.  We denote by
$\mathrm{Sing}(\mathcal{V})$ the singular locus of $\mathcal{V}$.

Let $d\geq 1$ be an integer. Define the 
\textit{algebra of fractional series} $\widetilde{\mathbf{C}\{ X\}}:= 
\displaystyle \mbox{lim}_{\stackrel{\longrightarrow}{N \geq 0}}
\mathbf{C}\{ X_{1}^{\frac{1}{N}},...,X_{d}^{\frac{1}{N}}\}$, where 
$X:=(X_{1},...,X_{d})$. 
If $m=(m_{1},...,m_{d})\in \mathbf{Q}_{+}^{d}$, we denote $X^{m}:=
X_{1}^{m_{1}}\cdots X_{d}^{m_{d}}$. 
If $\eta \in \widetilde{\mathbf{C}\{ X\}}$ can be written 
$\eta=X^{m}u(X)$, with $ m \in \mathbf{Q}_{+}^{d}$ and 
$u \in \widetilde{\mathbf{C}\{ X\}}, \: u(0,...,0)\neq 0$, we say that 
$\eta$ has a \textit{dominating exponent}.

\begin{defin} \label{qobj}
Let  $(\mathcal{S},0)$ be a germ of reduced equidimensional complex
space.  
The germ $(\mathcal{S},0)$ is called 
\textbf{quasi-ordinary} if there exists a finite morphism $\psi$ from 
$(\mathcal{S},0)$ to a smooth space of the same dimension, 
whose discriminant locus is contained in a hypersurface with 
normal crossings. Such a morphism $\psi$ is also called 
\textbf{quasi-ordinary}.
\end{defin}

For instance, all reduced germs of curves are quasi-ordinary with
respect to any finite 
morphism whose target is a smooth curve.

In the special case in which $\mathcal{S}$ is a $d$-dimensional
hypersurface germ,  
one can find local coordinates $X$ on the target space of $\psi$ 
such that the discriminant locus of $\psi$ is contained in
$\{X_{1}\cdots X_{d}=0\}$ and an element $Y$ in the maximal ideal 
of $\mathcal{O}_{\mathcal{S},0}$ such that $(\psi, Y)$ embeds
$(\mathcal{S},0)$ in  
$\mathbf{C}^{d}\times \mathbf{C}$. So $\psi$ appears as a map:
$$\psi:\mathcal{S} \rightarrow \mathbf{C}^{d},$$
which is \textit{unramified over} $(\mathbf{C}^{*})^{d}$. 
By the Weierstrass preparation 
theorem, the image of $\mathcal{S}$ by $(\psi, Y)$, 
identified in the sequel with 
$\mathcal{S}$, is defined by a unitary polynomial $f \in 
\mathbf{C}\{ X \}[Y]$. The discriminant locus of $\psi$ is defined 
by the discriminant $\Delta_{Y}(f)$ of $f$, which has therefore a
dominating exponent.

\begin{defin} \label{qopol}
Let $f \in \mathbf{C}\{ X \}[Y]$ be unitary. If $\Delta_{Y}(f)$ has a 
dominating exponent, we say that $f$ is \textbf{quasi-ordinary}. 
\end{defin}

The following theorem (see \cite{A 55},
\cite{L 88}), generalizes the theorem of Newton-Puiseux for plane
curves: 
 
\begin{theor} \label{JAbh} (Jung-Abhyankar) 
If $f \in \mathbf{C}\{ X\}[Y]$ is 
quasi-ordinary, then the set $R(f)$ of roots of $f$ embeds canonically 
in the algebra $\widetilde{\mathbf{C}\{ X\}}$.
\end{theor}
 
In the sequel, we consider $R(f)$ as a subset of 
$\widetilde{\mathbf{C}\{ X\}}$. Moreover, we suppose that $f$ is
\textit{irreducible}.  
Then all the differences of roots of $f$ have dominating exponents,
which are totally 
ordered for the componentwise order (see \cite{L 65}, 
\cite{L 88}). If $G$ is their number, denote them by 
$A_{1}<\cdots < A_{G}, \: A_{i}= (A^{1}_{i},...,A^{d}_{i}), \forall 
i \in \{1,...,G\}$. 

\begin{defin} \label{charexp}
We call the vectors $A_{1},...,A_{G} \in \mathbf{Q}_{+}^{d}$  
the \textbf{characteristic exponents} and the monomials
$X^{A_{1}},...,X^{A_{G}}$ the \textbf{characteristic monomials} of
$f$ or of $\psi$. 
\end{defin}

Some comments on the characteristic exponents follow the proposition
\ref{algfond}.

\section{Generalized Hirzebruch-Jung singularities} \label{genHJ}

In this section we recall some results about the normalization of
quasi-ordinary singularities and we define 
Hirzebruch-Jung singularities in any dimension. 

For details about toric geometry, see Oda \cite{O 88} and Fulton 
\cite{F 93}. 

We denote by $W_{0}=\mathbf{Z}^{d}$ the canonical $d$-dimensional
lattice, by $M_{0}=\mathbf{Z}^{d}$ its canonical dual and by
$\sigma_{0}$ the canonical regular cone of maximal dimension in $W_{0}$. 

Let $(\mathcal{S},0)$ be an irreducible $d$-dimensional quasi-ordinary
germ and let $\psi:(\mathcal{S},0)
\rightarrow (\mathbf{C}^{d},0)$ be a  finite morphism
unramified over $(\mathbf{C}^{*})^{d}$. We look at $\mathbf{C}^{d}$ as
the affine toric variety $\mathcal{Z}(W_{0}, \sigma_{0})$. Then the
fundamental group $\pi_{1}((\mathbf{C}^{*})^{d})$ can be canonically
identified with $W_{0}$ (see \cite{F 93}). 

Define:
$$W(\psi):=\psi_{*}\pi_{1}(\psi^{-1}((\mathbf{C}^{*})^{d})).$$
It is a subgroup of
$\pi_{1}((\mathbf{C}^{*})^{d})=W_{0}$. Moreover, $W(\psi)$
is of finite index in $W_{0}$, as $\psi$ is finite.  Consider the
affine toric variety 
$\mathcal{Z}(W(\psi), \sigma_{0})$ obtained by changing the lattice
from $W_{0}$ to $W(\psi)$. Denote by:
$$\gamma_{W_{0}:W(\psi)}:\mathcal{Z}(W(\psi), \sigma_{0}) \rightarrow 
    \mathcal{Z}(W_{0}, \sigma_{0})=\mathbf{C}^{d}$$
the canonical morphism associated to this change of lattice. We
proved topologically the following theorem in \cite{PP 01} and \cite{PP
  03'}. A more algebraic proof was given later by Aroca and Snoussi in
\cite{AS 03}. 

\begin{theor} \label{normtor}
One has
the following commutative diagram, in which $\nu$  
is a normalization morphism:
$$\xymatrix{
    (\mathcal{Z}(W(\psi), \sigma_{0}), 0)
         \ar[dr]_{\gamma_{W_{0}:W(\psi)}} 
             \ar[rr]^-{\nu} & 
       &  (\mathcal{S}, 0)
                 \ar[dl]^{\psi} \\
    & (\mathbf{C}^{d},0) & }$$
\end{theor}

In the special case in which $\mathcal{S}$ is a hypersurface germ, we
can express the lattice $W(\psi)$ using the characteristic 
exponents of $\psi$. In order to do this let us introduce, following
Lipman \cite{L 88}, the abelian groups
$M_{0}:= \mathbf{Z}^{d},\: M_{i}:= M_{i-1}+ \mathbf{Z}A_{i}, \: \forall 
i\in \{1,...,G\}$ and 
the successive indices $N_{i}:= (M_{i}:M_{i-1}), 
\: \forall 
i\in \{1,...,G\}$. Following
Gonz{\'a}lez P{\'e}rez \cite{GP 00} we consider also the dual 
lattices $W_{k}$ of the lattices $M_{k}$:
$$W_{k}:=\mathrm{Hom}(M_{k},\mathbf{Z}), \: \forall \: k \in \{1,...,G\}.$$
One has the inclusions:
$   M_{0} \subsetneq M_{1} \subsetneq \cdots \subsetneq M_{G}, \:
    W_{0} \supsetneq W_{1} \supsetneq \cdots \supsetneq W_{G}.$
The following proposition was proved in \cite{PP 01} and \cite{PP 03'}:

\begin{prop} \label{exprdiff}
Let $f \in \mathbf{C}\{X\}[Y]$ be an irreducible quasi-ordinary
polynomial and   
$\psi$ be the associated quasi-ordinary projection. Then  $W(\psi)= W_{G}.$
\end{prop}

Using this identification, theorem \ref{normtor} becomes:

\begin{corollary} \label{normtor1} (Gonz{\'a}lez P{\'e}rez) 
If $f$ is an irreducible quasi-ordinary polynomial defining the germ
$\mathcal{S}$, then 
one has the following commutative diagram, in which $\nu$ is a normalization  
morphism:
$$\xymatrix{
    (\mathcal{Z}(W_{G}, \sigma_{0}), 0)
         \ar[dr]_{\gamma_{W:W_{G}}} 
             \ar[rr]^-{\nu} & 
       &  (\mathcal{S}, 0)
                 \ar[dl]^{\psi} \\
    & (\mathbf{C}^{d},0) & }$$
\end{corollary}
\medskip

This theorem was first proved algebraically by Gonz{\'a}lez P{\'e}rez
in \cite{GP 00}, without passing through
proposition \ref{exprdiff}. It inspired our theorem \ref{normtor}.

Theorem \ref{normtor} and the fact that in dimension 2 Hirzebruch-Jung
singularities are precisely the normalizations of quasi-ordinary ones,
motivates us to introduce the following definition in arbitrary dimension:

\begin{defin} \label{defJH} 
The irreducible germ $(\mathcal{Z},0)$ of complex analytical space is 
called \textbf{a Hirzebruch-Jung singularity} if it is analytically 
isomorphic with the germ at the $0$-dimensional orbit of an 
affine toric variety defined by a maximal simplicial pair. 
\end{defin} 

One can give another definition of Hirzebruch-Jung singularities:

\begin{prop} \label{equidef}
Hirzebruch-Jung singularities are precisely the quasi-ordinary
singularities which are normal. 
\end{prop}

\textbf{Proof:} By theorem \ref{normtor}, each normal quasi-ordinary
singularity is a Hirze\-bruch-Jung one.

Conversely, let $(\mathcal{Z},0)$ be a Hirzebruch-Jung singularity,
according to definition \ref{defJH}. Then
$(\mathcal{Z},0)\simeq(\mathcal{Z}(\mathcal{W}, \sigma),0)$, where
$(\mathcal{W}, \sigma)$ is a maximal simplicial pair. Let
$v_{1},...,v_{d}$ be the primitive elements of $\mathcal{W}$ situated
on the edges of $\sigma$ and let $w_{1},...,w_{d}$ be a basis of the
lattice $\mathcal{W}$. The matrix transforming $(v_{1},...,v_{d})$ in
$(w_{1},...,w_{d})$ has rational coefficients. So, there is a number
$q \in \mathbf{N}^{*}$ such that the matrix transforming
$\Big(\dfrac{1}{q}v_{1},..., \dfrac{1}{q}v_{d}\Big)$ into
$(w_{1},...,w_{d})$ has 
integer coefficients. If $\mathcal{W}_{0}:=
\sum_{i=1}^{d}\mathbf{Z}\dfrac{1}{q}v_{i} \subset
\mathcal{W}_{\mathbf{Q}}$, then $\mathcal{W}$ is a sublattice of
finite index of $\mathcal{W}_{0}$ and $(\mathcal{W}_{0}, \sigma)$ is a
maximal regular pair. Denote also by $0$ the 0-dimensional orbit of
$\mathcal{Z}(\mathcal{W}_{0}, \sigma)$. Consider the toric morphism
obtained by changing the lattice:
$$\eta: \mathcal{Z}(\mathcal{W}, \sigma)\rightarrow
\mathcal{Z}(\mathcal{W}_{0}, \sigma)\simeq \mathbf{C}^{d}.$$
It is finite and unramified over the torus
$\mathcal{Z}(\mathcal{W}_{0}, \{0\})$. So, its germ \linebreak $\eta:
(\mathcal{Z}(\mathcal{W}, \sigma),0)\rightarrow  
(\mathcal{Z}(\mathcal{W}_{0}, \sigma),0)$ is quasi-ordinary. But 
$\mathcal{Z}(\mathcal{W}, \sigma)$ is normal, and the proposition is
proved.  \hfill $\Box$

\section{The analytical classification \\of Hirzebruch-Jung
  singularities} \label{mainres} 

In this section we prove our main result (theorem \ref{mthm}) which
classifies Hirzebruch-Jung singularities up to analytical isomorphism
by combinatorial data. It states that a maximal simplicial pair
$(\mathcal{W}, \sigma)$ can be reconstructed from the analytical type
of the germ $(\mathcal{Z}(\mathcal{W}, \sigma),0)$. Our essential tool
is \textit{the orbifold map} $\mu$ 
associated to $\mathcal{Z}$ (definition \ref{orbdef}). Then we 
give a normal form for 
maximal simplicial pairs $(\mathcal{W}, \sigma)$, once an ordering of
the edges of $\sigma$ was fixed (proposition \ref{normcone}). This
allows us to define \textit{the type} of a
Hirzebruch-Jung singularity (definition \ref{type}). 

Let $V$ be a complex finite dimensional vector space. An element of
$GL(V)$ is called a \textit{complex reflection} if its fixed-point set
is a hyperplane of $V$. A finite group $\Gamma \subset GL(V)$ is
called \textit{small} (see Prill \cite{P 67}) if it contains no complex
reflections.  

Let us recall a generalization of the Riemann
existence theorem (see \cite{BN 90}):

\begin{theor} \label{GRexist}(Grauert-Remmert)
Let $S$ be a connected normal complex space and $T\subset S$ a proper
closed 
analytical subset. Let $Y:=S-T$, let $X$ be a normal complex space and
$\phi:X \rightarrow Y$ be a ramified covering. Then $\phi$ extends to
a ramified covering $\widetilde{\phi}: \widetilde{X}\rightarrow S$ with
$\widetilde{X}$ normal if and only if the closure $\overline{B}$ in $S$
of the branch locus $B \subset Y$ of $\phi$ is an analytical subset in
$S$. In this case, the extension is unique.
\end{theor}

In the case it exists, we say that $\widetilde{\phi}$ is obtained by
\textit{Riemann extension} of $\phi$.

Let $(\mathcal{W},\sigma)$ be a maximal simplicial pair of dimension $d \geq
2$.  Denote by  $v_{i}, \: i \in \{1,...,d\}$ the primitive elements
of $\mathcal{W}$ situated on the edges of $\sigma$. Denote by 
  $\widetilde{\mathcal{W}}$ the sublattice of $\mathcal{W}$ generated by
  $v_{1},...,v_{d}$. Then $(\widetilde{\mathcal{W}}, \sigma)$ is a
  maximal regular pair. Consider the toric morphism: 
    $$\mu:
  \mathcal{Z}(\widetilde{\mathcal{W}}, \sigma) \rightarrow
  \mathcal{Z}(\mathcal{W}, \sigma),$$ 
  obtained by keeping the same cone $\sigma$ and by replacing the
  lattice $\mathcal{W}$ by $\widetilde{\mathcal{W}}$. In what follows, we
  will denote by 
  $(\mathcal{Z},0)$ and $(\widetilde{\mathcal{Z}},0)$ the complex analytic
  germs
  $(\mathcal{Z}(\mathcal{W}, \sigma), 0)$ and
  $(\mathcal{Z}(\widetilde{\mathcal{W}}, \sigma), 
  0)$, or sufficiently small representatives of them. Notice that
  $0=\mu^{-1}(0)$.  

\begin{prop} \label{applext}
The map $\mu$ is obtained by Riemann extension of the
universal covering map of the smooth part of $\mathcal{Z}(\mathcal{W},
\sigma)$.  
In particular, the restriction of $\mu$ over the germ
$(\mathcal{Z},0)$ depends only on the analytical structure of
$(\mathcal{Z},0)$.
\end{prop}

\textbf{Proof:} The proof of this proposition is also contained in the
section 6 of \cite{PP 03'}.

By general results of toric geometry (see \cite{O 88}, corollary
1.16), $\mu$ is the 
quotient map of $\mathcal{Z}(\widetilde{\mathcal{W}}, \sigma)$ by the
natural action of the finite group
$\mathcal{W}/\widetilde{\mathcal{W}}$. Moreover, in toric coordinates,
this action is linear, faithful, and does not contain complex
reflections. So, as a linear group
$\mathcal{W}/\widetilde{\mathcal{W}}$ is small. A rapid proof of this
fact will be given in the remark which follows the proof of theorem
\ref{mthm}.  

This shows that the locus $\mathrm{Fix}(\mu)$ of the fixed points
  of the elements of $\mathcal{W}/\widetilde{\mathcal{W}}$ distinct from
  the identity  has codimension at least 2 in
  $\mathcal{Z}(\widetilde{\mathcal{W}}, \sigma)$. Moreover,
  $\mu^{-1}(\mathrm{Sing}(\mathcal{Z}(\mathcal{W},\sigma)))\subset
  \mathrm{Fix}(\mu)$. As $\mathcal{Z}(\widetilde{\mathcal{W}},
  \sigma)$ is smooth, 
  the complement $\mathcal{Z}(\widetilde{\mathcal{W}},
  \sigma)-\mu^{-1}(\mathrm{Sing}(\mathcal{Z}(\mathcal{W},\sigma)))$ is
  simply connected, and so the restriction of $\mu$ over the smooth
  part of $\mathcal{Z}(\mathcal{W},\sigma)$ is a universal covering
  map. The uniqueness in theorem \ref{GRexist} implies the
  proposition. \hfill $\Box$ 

\medskip

Following a terminology used in \cite{DM 93}, we define:

\begin{defin}\label{orbdef}
 The morphism
$\mu$ obtained by Riemann extension of the universal covering map of
the smooth part of $(\mathcal{Z},0)$ is called \textbf{the orbifold
  map} associated to $(\mathcal{Z}, 0)$.  
\end{defin}

Denote by $\Gamma(\mathcal{Z})$ the \textit{group of covering
  transformations} of $\mu$ (in the terminology of \cite{D 88}), formed by
those analytical automorphisms $\phi: (\widetilde{\mathcal{Z}},0) \rightarrow
(\widetilde{\mathcal{Z}},0)$ which verify $\mu=\mu\circ \phi$. Consider
its action: 
\begin{equation} \label{repr}
 \Gamma(\mathcal{Z})\stackrel{\rho(\mathcal{Z})}{\longrightarrow}
 GL(\widetilde{m}/\widetilde{m}^{2})
\end{equation}
on the Zariski cotangent space of $\widetilde{\mathcal{Z}}$ at $0$. Here
$\widetilde{m}$ denotes the maximal ideal of $\widetilde{\mathcal{Z}}$ at
$0$. Being abelian, the group $\Gamma(\mathcal{Z})$ is canonically
isomorphic with the 
local fundamental group of $(\mathcal{Z},0)$. As an abstract
representation, $\rho(\mathcal{Z})$
is clearly determined by the analytical type of the germ $(\mathcal{Z},0)$. 
The previous analysis shows that the map (\ref{repr}) is a faithful 
$\mathbf{C}$-linear representation of 
$\Gamma(\mathcal{Z})$, whose image is small. 

More generally, consider a faithful finite-dimensional
$\mathbf{C}$-linear representation
  $$\Gamma \stackrel{\rho}{\longrightarrow} GL(V)$$
of a finite abelian group $\Gamma$, such that its image is
small. Denote $d=\mathrm{dim} \: V$. 
Choose a decomposition $V= 
E_{1}\oplus E_{2} \oplus \cdots \oplus E_{d}$ of $\rho$ as a sum of
irreducible (1-dimensional) representations. This is possible, since 
$\Gamma$ is
abelian (see \cite{S 67}). Denote by $\mathcal{E}$ this decomposition. 

For any $g\in \Gamma$ and any $k \in \{1,...,d\}$, $g$ acts on
$E_{k}$ by multiplication by a root of unity $e^{2i \pi w_{k}(g)}$. 
Here $w_{k}(g)\in \mathbf{Q}$ is well-defined modulo
$\mathbf{Z}$. Define then:
  $$w_{\mathcal{E}}(g):= (w_{1}(g),...,w_{d}(g))\in \mathbf{Q}^{d}.$$
This vector is well-defined modulo $\mathbf{Z}^{d}$. Define the
following over-lattice of $W_{0}=\mathbf{Z}^{d}$: 
\begin{equation} \label{ovlat}
 W_{\mathcal{E}}(\rho):= \mathbf{Z}^{d} + \Sigma_{g \in \Gamma}
 \mathbf{Z} w_{\mathcal{E}}(g).
\end{equation}
As the vectors $w_{\mathcal{E}}(g)$ are well-defined modulo
$\mathbf{Z}^{d}$, it is 
clear that $W_{\mathcal{E}}(\rho)$ does not depend on their
choices. Moreover, as the decomposition of a representation of a
finite group as direct
sum of irreducible ones is unique up to the order of the summands (see
\cite{S 67}), the pair $(W_{\mathcal{E}}(\rho),
\sigma_{0})$ is independent up to isomorphism of the choice of
decomposition $\mathcal{E}$. That is why we shall denote it shortly:
$$(W(\rho),\sigma_{0}).$$

Our main theorem is:

\begin{theor} \label{mthm}
The pairs $(\mathcal{W}, \sigma)$ and $(W(\rho(\mathcal{Z})), \sigma_{0})$ are
isomorphic.
\end{theor}

\textbf{Proof:} As a $\textbf{C}$-vector space,
$\textbf{C}[\widetilde{\mathcal{M}}\cap\check{\sigma}]$ is generated by the
monomials $\chi^{\widetilde{m}}$, with $\widetilde{m}\in
\widetilde{\mathcal{M}}\cap 
  \check{\sigma}$. The canonical action of $\mathcal{W}$ on these monomials is
  given by:
  \begin{equation} \label{act}
  (w, \chi^{\widetilde{m}})\rightarrow e^{2i\pi (w,
    \widetilde{m})}\chi^{\widetilde{m}}
  \end{equation}

Let $(\check{v_{1}},...,\check{v_{d}})$ be the basis of
$\widetilde{\mathcal{M}}$ dual to the basis $(v_{1},...,v_{d})$ of
$\widetilde{\mathcal{W}}$. Then the images of
$\chi^{\check{v_{1}}},...,\chi^{\check{v_{d}}}$ constitute a basis
of the $\mathbf{C}$-vector space $\widetilde{m}/\widetilde{m}^{2}$. For any $k
\in \{1,...,d\}$, denote by $F_{k}$ the subspace of
$\widetilde{m}/\widetilde{m}^{2}$ generated by the image of
$\chi^{\check{v_{k}}}$. Denote by $\mathcal{F}$ the 
decomposition $\widetilde{m}/\widetilde{m}^{2}= 
F_{1}\oplus F_{2} \oplus \cdots \oplus F_{d}$. Then, by formula
(\ref{act}), for any $w\in  
\mathcal{W}$ and any $k\in \{1,...,d\}$, $w$ acts on $F_{k}$ by multiplication
with $e^{2i\pi (w, \check{v_{k}})}$. If $g(w)$ denotes the image of
$w$ in the group $\Gamma(\mathcal{Z})\simeq
\mathcal{W}/\widetilde{\mathcal{W}}$, this shows that: 
$$w_{\mathcal{F}}(g(w))=((w, \check{v_{1}}),...,(w, \check{v_{d}})),$$
and so, by formula (\ref{ovlat}): 
 $$(W(\rho(\mathcal{Z})), \sigma_{0})\simeq (\mathbf{Z}^{d}+\sum_{w\in
   \mathcal{W}}\mathbf{Z}((w, \check{v_{1}}),...,(w, \check{v_{d}})),
 \sigma_{0}).$$
If we express the pair $(\mathcal{W}, \sigma)$ using the basis
$(v_{1},...,v_{d})$ of the associated $\mathbf{Q}$-vector space, we
get the isomorphism:
$$\begin{array}{ll}
(\mathcal{W}, \sigma)& =(\sum_{k=1}^{d}\mathbf{Z}v_{k}+\sum_{w\in
  \mathcal{W}}\mathbf{Z}((w, \check{v_{1}})v_{1}+\cdots+ (w,
\check{v_{d}})v_{d}), \sum_{k=1}^{d}\mathbf{R}_{+}v_{k})\\
  & \simeq (\mathbf{Z}^{d}+\sum_{w\in
  \mathcal{W}}\mathbf{Z}((w, \check{v_{1}}),...,(w, \check{v_{d}})),
\sigma_{0})
  \end{array}$$
which proves the theorem. \hfill $\Box$

\medskip

\noindent \textbf{Remark: } The constructions done in the previous
proof  show easily that 
the image of the group $\Gamma(\mathcal{Z})\simeq
\mathcal{W}/\widetilde{\mathcal{W}}$ by the representation
$\rho(\mathcal{Z})$ is small. Suppose this is false and consider $w
\in \mathcal{W}$ such that $g(w)$ acts on
$\widetilde{m}/\widetilde{m}^{2}$ as a complex reflection. Consider
again the basis of $\widetilde{m}/\widetilde{m}^{2}$ formed by the images of
$\chi^{\check{v_{1}}},...,\chi^{\check{v_{d}}}$. Possibly after
reordering it, we can
suppose that $(w, \check{v}_{i})\in \mathbf{Z},\: \forall
\: i \in \{1,...,d-1\}$ and $(w, \check{v}_{d})\notin \mathbf{Z}$. As
$w=\sum_{i=1}^{d} (w, \check{v}_{i})v_{i}$, this implies that $(w,
\check{v}_{d})v_{d}\in \mathcal{W}$. As $(w, \check{v}_{d})\notin
\mathbf{Z}$, this contradicts the fact that $v_{d}$ is a primitive
element of $\mathcal{W}$. 
\medskip

The following proposition shows that a representation $\rho$ and the
pair \linebreak $(W(\rho), \sigma_{0})$ it determines contain equivalent
information. 

\begin{prop} \label{repcomb}
Let $\Gamma \stackrel{\rho}{\longrightarrow} GL(V)$ be faithful
finite-dimensional $\mathbf{C}$-linear representation 
of a finite abelian group $\Gamma$ whose image is
small. If $(\mathcal{Z},0)$ denotes the Hirzebruch-Jung singularity
defined by $(W(\rho), \sigma_{0})$, then the representations $\rho$
and $\rho(\mathcal{Z})$ are isomorphic. 
\end{prop}

\textbf{Proof: } Choose $\mathcal{E}$, an arbitrary decomposition of
$\rho$ as a sum of irreducible representations. As $\rho(\Gamma)$ is
small, we see that no $d$-tuple $w_{\mathcal{E}}(g)$, with $g\neq 1$,
is contained on a line defined by an edge of $\sigma_{0}$. This shows
that the $d$-tuples modulo $\mathbf{Z}^{d}$ can be recovered from
$(W(\rho), \sigma_{0})$, simply by expressing the elements of
$W(\rho)$ in terms of the primitive elements situated on the edges of
$\sigma_{0}$. Moreover, there is a bijection between the elements of
$\Gamma$ and the set of these tuples in
$\mathbf{Q}^{d}/\mathbf{Z}^{d}$. Associate then to
$w_{\mathcal{E}}(g)$ the vector $(e^{2i \pi w_{1}(g)},..., e^{2i \pi
  w_{d}(g)})\in (\mathbf{C}^{*})^{d}$. This map is injective and
invariant modulo $\mathbf{Z}^{d}$. We get immediately the
proposition. \hfill $\Box$
\medskip

As $(W(\rho(\mathcal{Z})), \sigma_{0})$ is determined by the
analytical type of $(\mathcal{Z},0)$, an immediate corollary of the
theorem \ref{mthm} is the announced analytical classification of
Hirzebruch-Jung singularities:

\begin{corollary} \label{isom}
Let $\mathcal{Z}$ and $\mathcal{Z'}$ be two toric
varieties defined by maximal simplicial pairs. Denote by $0$ and $0'$
their closed orbits. Then the Hirzebruch-Jung singularities 
$(\mathcal{Z},0)$ and $(\mathcal{Z}',0')$ are isomorphic as germs of
complex analytical varieties if and only if $\mathcal{Z}$ and
$\mathcal{Z'}$ are isomorphic as toric varieties. 
\end{corollary}

The theorem \ref{mthm} and its corollary show that in order to
describe the analytical type of a Hirzebruch-Jung singularity, it is
enough to describe the combinatorial type of the pair $(\mathcal{W}, \sigma)$
associated to it. In the following proposition we give a normal form
for such a pair, once an ordering of the edges of $\sigma$ is 
fixed. We will denote by "$\prec$" such an ordering.

\begin{prop} \label{normcone}
Let $(\mathcal{W}, \sigma)$ be a maximal simplicial pair of dimension $d$.  Let
$v_{1},...,v_{d}$ be the primitive elements of $\mathcal{W}$ situated on the
  edges of $\sigma$, once an ordering $\prec$ of them is chosen.
  Then, 
there exists a unique basis $(e_{1},...,e_{d})$ of $\mathcal{W}$ such that the
vectors $(v_{1},...,v_{d})$ can be written as:
\begin{equation} \label{rel1}
  \left\{ \begin{array}{l}
     v_{1}=e_{1}\\
     v_{2}=-\alpha_{1,2} e_{1} +\alpha_{2,2} e_{2}\\
     v_{3}=-\alpha_{1,3} e_{1}-\alpha_{2,3} e_{2}+\alpha_{3,3} e_{3}\\
     ...................................................\\
     v_{d}=-\alpha_{1,d}e_{1}-\cdots-\alpha_{d-1,d}
     e_{d-1}+\alpha_{d,d} e_{d} 
          \end{array} \right.
\end{equation}
with $0\leq \alpha_{i,j}<\alpha_{j,j}$ for all $1\leq i<j\leq d$. 
\end{prop}

\textbf{Proof:}  If the given relations are verified, then $\forall k
\in \{1,...,d\}$, the vectors $e_{1},...,e_{k}$ are elements of the
lattice $\mathcal{W}\cap(\sum_{i=1}^{k}\mathbf{Q}v_{i})$. Moreover,
they form a basis of it, as $(e_{1},...,e_{d})$ is a basis of
$\mathcal{W}$. So, in order to prove the existence and the unicity of
$(e_{1},..., e_{d})$ once the conditions $0\leq
\alpha_{i,j}<\alpha_{j,j}$ are imposed, we will restrict to $d$-tuples
of vectors such 
that $(e_{1},...,e_{k})$ is a basis of \linebreak 
$\mathcal{W}\cap(\sum_{i=1}^{k}\mathbf{Q}v_{i}), \: \forall \: k \in
\{1,...,d\}$.

It is clear that $e_{1}$ exists and is unique verifying the first
relation. 

Suppose that $(e_{1},...,e_{k-1})$ is a basis of the rank $(k-1)$
lattice $\mathcal{W}\cap(\sum_{i=1}^{k-1}\mathbf{Q}v_{i})$ that
verifies the first $(k-1)$ relations of (\ref{rel1}), where $k \geq 2$. Choose
$\widetilde{e}_{k}\in \mathcal{W}$ such that $(e_{1},...,e_{k-1},
\widetilde{e}_{k})$ is a basis of the rank $k$ lattice
$\mathcal{W}\cap(\sum_{i=1}^{k}\mathbf{Q}v_{i})$. This is possible, as
the quotient
$(\mathcal{W}\cap(\sum_{i=1}^{k}\mathbf{Q}v_{i}))/
(\mathcal{W}\cap(\sum_{i=1}^{k-1}\mathbf{Q}v_{i}))$ has no torsion. 
Then one can write:
$$v_{k}=-\widetilde{\alpha}_{1,k}e_{1}- \cdots
-\widetilde{\alpha}_{k-1,k}e_{k-1}+\widetilde{\alpha}_{k,k}\widetilde{e}_{k}$$ 
with $\widetilde{\alpha}_{i,k}\in \mathbf{Z}, \: \forall \: i \in
\{1,...,k\}$. 

If $e_{k}'$ is such that $(e_{1},...,e_{k-1}, e_{k}')$ is also a basis
of $\mathcal{W}\cap(\sum_{i=1}^{k}\mathbf{Q}v_{i})$, then:
$$\widetilde{e}_{k}=\epsilon
e_{k}'+\lambda_{k-1}e_{k-1}+\cdots+\lambda_{1}e_{1},$$  
where $\epsilon \in \{+1, -1\}$ and $\lambda_{i}\in \mathbf{Z}, \: \forall
\: i \in \{1,...,k-1\}$. So:
$$v_{k}= - (\widetilde{\alpha}_{1,k}-\lambda_{1}\widetilde{\alpha}_{k,k})e_{1}
-\cdots -
(\widetilde{\alpha}_{k-1,k}-\lambda_{k-1}\widetilde{\alpha}_{k,k})e_{k-1}+ 
\epsilon \widetilde{\alpha}_{k,k}e_{k}'.$$

The number $\epsilon$ is uniquely determined by the condition
$\epsilon \widetilde{\alpha}_{k,k} >0$. Then, $\forall \: i \in
\{1,...,k-1\}$, the integer $\lambda_{i}$ is clearly uniquely determined by
the condition $0 \leq
\widetilde{\alpha}_{i,k}-\lambda_{i}\widetilde{\alpha}_{k,k}<\epsilon
\widetilde{\alpha}_{k,k}$. 

So, there exists a unique $e_{k}$ such that $(e_{1},...,e_{k-1},
e_{k})$ verify the first $k$ relations of (\ref{rel1}). This proves
the proposition by induction. \hfill $\Box$

\medskip

We denote by $\mathcal{B}(\mathcal{W}, \sigma, \prec)$ the basis
$(e_{1},...,e_{d})$ of $W$ and by $\mathfrak{m}(\mathcal{W}, \sigma,
\prec)$ the matrix:
$$\left( \begin{array}{llll}
           1 & -\alpha_{1,2} & \cdots & -\alpha_{1,d}\\
           0 & \alpha_{2,2}  & \cdots & -\alpha_{2,d}\\
           \vdots &  \vdots & \ddots & \vdots \\
           0 & 0 & \cdots & \alpha_{d,d}
         \end{array} \right)$$

\begin{defin} \label{type}
Let $(\mathcal{Z},0)$ be a Hirzebruch-Jung singularity isomorphic with
$(\mathcal{Z}(\mathcal{W}, \sigma),0)$, where $(\mathcal{W}, \sigma)$
is a maximal 
simplicial pair. If $\prec$ is an ordering of the edges of $\sigma$,
we say that $(\mathcal{Z},0)$ \textbf{is of type} 
$\mathfrak{m}(\mathcal{W}, \sigma, \prec)$. 
\end{defin}

We see that there is a finite ambiguity in the definition of the type
of $(\mathcal{Z},0)$. Indeed, there are $d!$ possible orderings, and
so $d!$ possible matrices $\mathfrak{m}(\mathcal{W}, \sigma, \prec)$.
\medskip

\noindent \textbf{Remark: } It would be interesting to find a method
to decide if two matrices 
correspond to the same pair $(\mathcal{W}, \sigma)$ but to distinct choices of
the ordering of the edges of $\sigma$. Such a method is known
classically in dimension 2 (see 
proposition \ref{amb}).

\section{A normalization algorithm for quasi-ordinary hypersurface
  singularities }
 \label{sectalg}

We proved in \cite{PP 01} (see also 
\cite{PP 03}) an algorithm for computing the Hirzebruch-Jung type of
the normalization of a quasi-ordinary singularity of hypersurface in
$\mathbf{C}^{3}$. In this section we generalize it to arbitrary
dimensions (proposition \ref{algfond}). In order to do it, we need to
give a normal form for sublattices $\overline{\mathcal{W}}$ of finite
index of a lattice $\mathcal{W}$, once a basis of $\mathcal{W}$ is
fixed (proposition \ref{normlat}). As an important intermediate
result, we give an algorithm of computation of this normal form when
$\overline{\mathcal{W}}$ is defined by a congruence (lemma
\ref{cal3}). 

By corollary \ref{normtor1}, the normalization of the germ
$(\mathcal{S},0)$ defined by an irreducible quasi-ordinary polynomial
is a Hirzebruch-Jung singularity of type $\mathfrak{m}(W_{G},
\sigma_{0}, \prec_{0})$. This shows that we need to compute the
sublattice $W_{G}$ of $W_{0}$. We will first prove a proposition similar to
proposition \ref{normcone}, which gives a normal form to a sublattice
$\overline{\mathcal{W}}$ of finite index of a given lattice
$\mathcal{W}$, once a basis of $\mathcal{W}$ has been fixed.

\begin{prop} \label{normlat}
 Let $(\mathcal{W},\sigma)$ be a maximal regular pair of dimension
 $d$. Let 
 $(w_{1},...,w_{d})$ be the primitive elements of $\mathcal{W}$
 situated on the edges of $\sigma$, once an ordering $\prec$ of them has been
 chosen. If $\overline{\mathcal{W}}$ is a sublattice of finite
 index of $\mathcal{W}$, then there exists a unique basis
 $(\overline{w}_{1},...,\overline{w}_{d})$ of $\overline{\mathcal{W}}$
 such that: 
\begin{equation} \label{rel2}
  \left\{ \begin{array}{l}
     \overline{w}_{1}=r_{1,1}w_{1}\\
     \overline{w}_{2}=r_{1,2} w_{1}+r_{2,2} w_{2}\\
     \overline{w}_{3}=r_{1,3} w_{1}+r_{2,3} w_{2}+r_{3,3} w_{3}\\
     ...............................................\\
     \overline{w}_{d}=r_{1,d} w_{1}+r_{2,d} w_{2}+\cdots+r_{d,d}w_{d} 
        \end{array} \right.
\end{equation}
with $0\leq r_{i,j} < r_{i,i}$ for all $1 \leq i < j \leq d$. 
\end{prop}

\textbf{Proof:} If the given relations are verified, then $\forall \:
k \in \{1,...,d\}$, the vectors
$\overline{w}_{1},...,\overline{w}_{k}$ are elements of the lattice
$\overline{\mathcal{W}}\cap (\sum_{i=1}^{k}\mathbf{Z}
w_{i})$. Moreover, they form a basis of it, as
$(\overline{w}_{1},...,\overline{w}_{d})$ is a basis of
$\overline{\mathcal{W}}$. So, in order to prove the existence and the
unicity of $(\overline{w}_{1},...,\overline{w}_{d})$ once the
conditions $0\leq r_{i,j} < r_{i,i}$ are imposed, we will restrict
to $d$-tuples of vectors such that
$(\overline{w}_{1},...,\overline{w}_{k})$ is a basis of \linebreak 
$\overline{\mathcal{W}}\cap (\sum_{i=1}^{k}\mathbf{Z} w_{i}), \:
\forall \: k \in \{1,...,d\}$. 

Consider the rank 1 lattice
$\overline{\mathcal{W}}\cap \mathbf{Z} w_{1}$. It has a unique
generator of the form $r_{1,1}w_{1}$, with $r_{1,1}>0$. Set
$\overline{w}_{1}:= r_{1,1}w_{1}$. 

Suppose now that $(\overline{w}_{1},...,\overline{w}_{k-1})$ is a
basis of the rank $(k-1)$ lattice \linebreak $\overline{\mathcal{W}}\cap
(\sum_{i=1}^{k-1}\mathbf{Z} w_{i})$ that verifies the $(k-1)$ first
relations of (\ref{rel2}), where $k \geq 2$. Choose
$\widetilde{w}_{k}$ such that 
$(\overline{w}_{1},...,\overline{w}_{k-1}, \widetilde{w}_{k})$ is a 
basis of the rank $k$ lattice $\overline{\mathcal{W}}\cap
(\sum_{i=1}^{k}\mathbf{Z} w_{i})$. Then one can write:
$$\widetilde{w}_{k}=\widetilde{r}_{1,k}w_{1}+\cdots+\widetilde{r}_{k,k}w_{k},$$
with $\widetilde{r}_{i,k}\in \mathbf{Z}, \: \forall \: i \in
\{1,...,k\}$. 

If $w_{k}'$ is such that $(\overline{w}_{1},...,\overline{w}_{k-1},
w_{k}')$ is also a basis of $\overline{\mathcal{W}}\cap
(\sum_{i=1}^{k}\mathbf{Z} w_{i})$, then:
$$w_{k}'= \epsilon \widetilde{w}_{k} +\lambda_{k-1}\overline{w}_{k-1}+ \cdots +
\lambda_{1}\overline{w}_{1}$$ 
where $\epsilon \in \{+1, -1\}$ and $\lambda_{i}\in \mathbf{Z}, \: \forall
\: i \in \{1,...,k-1\}$. So:
$$\begin{array}{ll}
    w_{k}' & = \epsilon \sum_{i=1}^{k}\widetilde{r}_{i,k} w_{i} +
    \sum_{j=1}^{k-1}\lambda_{j}( \sum_{i=1}^{j}r_{i,j}w_{i})=\\
        & = \sum_{i=1}^{k-1}(\epsilon \widetilde{r}_{i,k}+
        \sum_{j=i}^{k-1}\lambda_{j}r_{i,j}) w_{i} + \epsilon
        \widetilde{r}_{k,k} 
        w_{k}
   \end{array}$$

The number $\epsilon$ is uniquely determined by the condition
$\epsilon \widetilde{r}_{k,k}>0$. Then, $\lambda_{k-1}$ is uniquely determined
by the condition $0 \leq \epsilon
\widetilde{r}_{k-1,k}+\lambda_{k-1}r_{k-1,k-1} 
<r_{k-1,k-1}$, so $\lambda_{k-2}$ is uniquely determined by the condition 
$0 \leq \epsilon
\widetilde{r}_{k-2,k}+\lambda_{k-2}r_{k-2,k-2}+\lambda_{k-1}r_{k-2,k-1} 
<r_{k-2,k-2}$. Keeping like this, we see that $(\epsilon, \lambda_{k-1},
\lambda_{k-2},...,\lambda_{1})$ are uniquely determined by the conditions
$\epsilon \widetilde{r}_{k,k}>0$ and $0 \leq \epsilon \widetilde{r}_{i,k}+
        \sum_{j=i}^{k-1}\lambda_{j}r_{i,j} < r_{i,i}, \: \forall \: i \in
        \{1,...,k-1\}$. 

 This proves the proposition by
induction. \hfill $\Box$ 

\medskip

We denote by $\mathfrak{B}(\mathcal{W}, \sigma, \prec;
\overline{\mathcal{W}})$ the basis
$(\overline{w}_{1},...,\overline{w}_{d})$ of $\overline{\mathcal{W}}$
and by \linebreak $\mathfrak{m}(\mathcal{W}, \sigma, \prec;
\overline{\mathcal{W}})$ the matrix:
$$\left( \begin{array}{llll}
          r_{1,1} & r_{1,2} & \cdots & r_{1,d}\\
          0 & r_{2,2} & \cdots & r_{2,d}\\
          \vdots & \vdots & \ddots & \vdots\\
          0 & 0& \cdots & r_{d,d}
         \end{array} \right).$$

Suppose now that the relations (\ref{rel2}) are verified but perhaps
without satisfying the conditions $0\leq r_{i,j} < r_{i,i}$. Let
$\mathfrak{m}$ be the matrix $(r_{i,j})_{i,j}$. Denote by:
$$n_{s}(\mathfrak{m})$$
the matrix $\mathfrak{m}(\mathcal{W}, \sigma, \prec;
\overline{\mathcal{W}})$. Here "s" is the initial letter of
"sublattice". This alludes to the fact that one has to choose the
base of the sublattice $\overline{\mathcal{W}}$ of $\mathcal{W}$. The
proof of proposition \ref{normlat} gives 
an algorithm of computation of $n_{s}(\mathfrak{m})$ starting
from the knowledge of $\mathfrak{m}$. 

Analogously, if the relations (\ref{rel1}) are verified but perhaps
without satisfying the conditions $0 \leq \alpha_{i,j} <\alpha_{j,j}$,
and $\mathfrak{m}$ denotes the matrix transforming $(e_{1},...,e_{d})$
into $(v_{1},...,v_{d})$, we denote by:
$$n_{a}(\mathfrak{m})$$
the matrix $\mathfrak{m}(\mathcal{W}, \sigma, \prec)$. Here "a" is the
initial letter of "ambient lattice", it alludes to the fact that one
has to choose the base of the ambient lattice $\mathcal{W}$. The
proof of proposition \ref{normcone} gives 
an algorithm of computation of $n_{a}(\mathfrak{m})$ starting
from the knowledge of $\mathfrak{m}$.

 If $t\in
\mathbf{Q}$, one can write in a unique way $t=\frac{p}{q}$ with
$\mathrm{gcd}(p,q)=1$ and $q >0$. Define the \textit{numerator} and
the \textit{denominator} of $t$ by:
$$\mathrm{num}(t):=p,$$
 $$ \mathrm{den}(t):=q.$$

The following lemma relates the
normal forms of the propositions \ref{normcone} and \ref{normlat}:

\begin{lemm} \label{transf}
If $(\mathcal{W}, \sigma)$ is a maximal simplicial pair, that $\prec$
is an ordering of the edges of $\sigma$ and that
$\overline{\mathcal{W}}$ is a sublattice of finite index of
$\mathcal{W}$, then the matrix $\mathfrak{m}(\mathcal{W}, \sigma, \prec;
\overline{\mathcal{W}})$ determines the matrix
$\mathfrak{m}(\overline{\mathcal{W}}, \sigma, \prec)$. 
\end{lemm}

\textbf{Proof: } As $\mathfrak{m}(\mathcal{W}, \sigma, \prec;
\overline{\mathcal{W}})$ is upper triangular, so is its inverse. But
unlike the entries of $\mathfrak{m}(\mathcal{W}, \sigma, \prec;
\overline{\mathcal{W}})$, the entries $t_{i,j}$ of
$\mathfrak{m}(\mathcal{W}, \sigma, \prec;
\overline{\mathcal{W}})^{-1}$ are not in general integers. As
$w_{j}=\sum_{i=1}^{j}t_{i,j}\overline{w}_{i}, \: \forall \: j \in
\{1,...,d\}$, one sees that $(d_{j}w_{j})_{1 \leq j \leq d}$ are the
primitive elements of $\overline{\mathcal{W}}$ situated  on the edges
of $\sigma$, their ordering $\prec$ being the same as before. We have
denoted:
$$d_{j}:= \mathrm{lcm}(\mathrm{den}(t_{1,j}),...,
\mathrm{den}(t_{j,j})), \: \forall \: j \in \{1,...,d\}.$$
This shows that:
\begin{equation} \label{normtrans}
  \mathfrak{m}(\overline{\mathcal{W}}, \sigma,
  \prec)=n_{a}((d_{j}t_{i,j})_{i,j}).
\end{equation}
$\:$ \hfill $\Box$

\medskip

Let $\mathcal{M}, \mathcal{W}$ be two dual rank $d$ lattices endowed
with dual basis $(u_{1},...,u_{d})$, respectively
$(w_{1},...,w_{d})$. Denote by $\sigma$ the cone spanned by
$(w_{1},...,w_{d})$. Consider $a \in \mathcal{M}_{\mathbf{Q}}$ and let
$\mathcal{W}(a)$ be the sublattice of $\mathcal{W}$ dual to
$\mathcal{M}(a):= \mathcal{M} + \mathbf{Z}a$, i.e. 
$\mathcal{W}(a):= \mathrm{Hom}(\mathcal{M} + \mathbf{Z}a,
\mathbf{Z}).$ 
 
If we write $a=\sum_{i=1}^{d}a^{i}u_{i}$, with $a^{1},...,a^{d}
\in \mathbf{Q}$, then:
\begin{equation} \label{cal1}
  \mathcal{W}(a)= \{ w \in \mathcal{W},\: (w,a) \in \mathbf{Z}\}= \{
  \sum_{i=1}^{d}c_{i}w_{i}, \: \sum_{i=1}^{d}c_{i}a^{i} \in
  \mathbf{Z}\}.
\end{equation}

\noindent \textbf{Remark: } The relation $\sum_{i=1}^{d}c_{i}a^{i} \in
  \mathbf{Z}$ can also be written $\sum_{i=1}^{d}(l_{d}a^{i})c_{i}
  \equiv 0 \:(\mathrm{mod}\: l_{d})$, where $l_{d}:=
  \mathrm{lcm}(\mathrm{den}(a^{1}), ..., \mathrm{den}(a^{d}))$. So,
  $\mathcal{W}(a)$ can be seen as a sublattice of finite index of
  $\mathcal{W}$ \textit{defined by a congruence}. 
\medskip

In the sequel we will also denote:
$$\mathfrak{m}(a^{1},...,a^{d}):= \mathfrak{m}(\mathcal{W}, \sigma, \prec;
\mathcal{W}(a)).$$

The following lemma gives an algorithm which computes the matrix \linebreak
$\mathfrak{m}(a^{1},...,a^{d})$ starting from the values of
$a^{1},...,a^{d}$.

\begin{lemm} \label{cal3}
Consider the matrix $\mathfrak{m}(a^{1},...,a^{d})=
(r_{i,j})_{i,j}$ and introduce the numbers 
$l_{k}:=\mathrm{lcm}(\mathrm{den}(a^{1}),...,\mathrm{den}(a^{k})), \:
\forall \: k \in \{1,...,d\}$. 
Then:
$$   r_{k,k} = \dfrac{l_{k}}{l_{k-1}}, \: \forall \: k \in
  \{1,...,d\}.$$
Moreover, for any $k \in \{1,...,d\}$ and any $j \in \{1,...,k-1\}$,
one has the equivalent relations:
$$\begin{array}{l}
     \sum_{i=j}^{k}l_{j}a^{i}r_{i,k}\equiv 0 \: (\mathrm{mod} \:
     r_{j,j})\\
       r_{j,k}=\left\{ \begin{array}{ll}
              -( \sum_{i=j+1}^{k}l_{j}a^{i}r_{i,k})(l_{j}a^{j})^{-1}
              \: \:        
                          \mbox{in} \:
                          \mathbf{Z}/r_{j,j}\mathbf{Z},  &
                          \mbox{ if } r_{j,j} \neq 1 \\
                        0, & \mbox{ if } r_{j,j} = 1
                       \end{array} \right.
  \end{array}$$
\end{lemm}

\textbf{Proof:} Denote by $t_{k}(a):=\sum_{j=1}^{k}a^{j}u_{j}$ the
$k$-truncation of $a$, for all $k \in \{1,...,d\}$. 

One knows (see the proof of lemma \ref{normlat}) that
$(\overline{w}_{1},...,\overline{w}_{k})$ is a basis of \linebreak 
$\mathcal{W}(a) \cap
(\sum_{j=1}^{k}\mathbf{Z}w_{j})=\mathcal{W}(t_{k}(a))$. This shows
that  $\prod_{j=1}^{k}r_{j,j}=(\mathcal{W}:
\mathcal{W}(t_{k}(a)))$. But, as the pairs of lattices $\mathcal{M},
\mathcal{W}$ and $\mathcal{M}(t_{k}(a)),\mathcal{W}(t_{k}(a))$ are
in duality, one has the equality of indices: $( \mathcal{W}:
\mathcal{W}(t_{k}(a)))= (\mathcal{M}(t_{k}(a)):\mathcal{M})$. This
last index is equal to the order of $t_{k}(a)$ in
$\mathcal{M}(t_{k}(a))/ \mathcal{M}$, which is obviously equal to
$l_{k}$. This implies:
$$\prod_{j=1}^{k}r_{j,j}=l_{k}$$    
which proves the first equalities.

Let us fix now $k \in \{2,...,d\}$ and $j \in \{1,...,k\}$. The
relation (\ref{cal1}) implies $\sum_{i=1}^{k}a^{i}r_{i,k}\in
\mathbf{Z}$. Multiplying this relation by $l_{j}$, we get: 
\begin{equation} \label{rel5}
  (\sum_{i=1}^{j-1}l_{j}a^{i}r_{i,k})+(\sum_{i=j}^{k}l_{j}a^{i}r_{i,k})\in
  l_{j}\mathbf{Z} \subset r_{j,j}\mathbf{Z}.
\end{equation}

But $\forall i\in \{1,...,j-1\}, \: l_{j}a^{i}=
r_{j,j}(l_{j-1}a^{i})\in r_{j,j}\mathbf{Z}$, as
$l_{j-1}a^{i} \in \mathbf{Z}$ by the definition of $l_{j-1}$. This
shows that $ \sum_{i=1}^{j-1}l_{j}a^{i}r_{i,k}\in r_{j,j}\mathbf{Z}$,
and (\ref{rel5}) implies:
\begin{equation} \label{rel6}
  \sum_{i=j}^{k}l_{j}a^{i}r_{i,k}\in
   r_{j,j}\mathbf{Z}.
\end{equation}
This is one of the forms in which were written the second relations of
the lemma.

Formula (\ref{rel6}) can also be written:
\begin{equation} \label{rel7}
  l_{j}a^{j}r_{j,k} +\sum_{i=j+1}^{k}l_{j}a^{i}r_{i,k}\in
   r_{j,j}\mathbf{Z}.
\end{equation}
As $\sum_{i=j+1}^{k}l_{j}a^{i}r_{i,k}=
\dfrac{\sum_{i=j+1}^{k}l_{j+1}a^{i}r_{i,k}}{r_{j+1,j+1}}$, relation
(\ref{rel6}) at the order $j+1$ implies that
$\sum_{i=j+1}^{k}l_{j}a^{i}r_{i,k} \in \mathbf{Z}$. Moreover,
$\mathrm{gcd}(l_{j}a^{j}, r_{j,j})=1$. Indeed, if $p$ is a prime
number dividing $r_{j,j}=\dfrac{l_{j}}{l_{j-1}}$, then $p \mid
\mathrm{den}(a^{j})$ and $p \nmid \dfrac{l_{j}}{\mathrm{den}(a^{j})}$.
As $\mathrm{gcd}(\mathrm{den}(a^{j}), \mathrm{num}(a^{j}))=1$, we also
have $p \nmid \mathrm{num}(a^{j})$, and so $p \nmid (l_{j}a^{j})$.
This 
shows that $l_{j}a^{j}$ is invertible in the ring
$\mathbf{Z}/r_{j,j}\mathbf{Z}$ if $r_{j,j}\neq 1$, and  from relation
(\ref{rel7}) we get the last formulae of 
the lemma. If $r_{j,j}=1$, as $0 \leq r_{j,k} <r_{j,j}$ we get
$r_{j,k}=0$. \hfill $\Box$

\medskip

\noindent \textbf{Remark: } The previous lemma shows that once
$r_{1,1}, ...,r_{d,d}$ are 
computed, one has to compute the entries of the $k$-th column in the
order: $r_{k-1,k}, r_{k-2,k},...,r_{1,k}$. As $0 \leq
r_{i,j}<r_{i,i}\: \forall \: 1 \leq i <j \leq d$, the entries
$r_{i,j}$ of the matrix $\mathfrak{m}(a^{1},...,a^{d})$ are completely
determined by the congruences of the lemma.

\medskip

Suppose now that $g \geq 1$ and $a_{1},...,a_{g}$ is a sequence of vectors of
$\mathcal{M}_{\mathbf{Q}}$. Define for all $k \in \{1,...,g\}$:
$$\begin{array}{l} 
  \mathcal{M}_{k}:= \mathcal{M}+\mathbf{Z}a_{1}+\cdots+
  \mathbf{Z}a_{k}\\
    \mathcal{W}_{k}:= \mathrm{Hom}(\mathcal{M}_{k}, \mathbf{Z}).
  \end{array}$$

We denote  $(r_{i,j}^{k})_{i,j}=\mathfrak{m}( \mathcal{W}, \sigma, \prec;
\mathcal{W}_{k})$. Write $a_{k}=a_{k}^{1}u_{1}+\cdots +
a_{k}^{d}u_{d}$, with $a_{k}^{1},...,a_{k}^{d}\in
\mathbf{Q}$. Introduce also the basis  
$\mathfrak{B}_{k}=(w_{1}^{k},...,w_{d}^{k}):=\mathfrak{B}( \mathcal{W},
\sigma, \prec; 
\mathcal{W}_{k})$. Denote by $\prec_{k}$ its ordering, deduced
canonically from $\prec$, and by $\sigma_{k}$ the cone generated by
$\mathfrak{B}_{k}$. Then:

$$\begin{array}{ll}
    \mathcal{W}_{k} & = \{ w\in \mathcal{W}_{k-1}, \: (w, a_{k})\in
    \mathbf{Z}\}=\\
        & =\{ w= \sum_{i=1}^{d}c_{i}w_{i}^{k-1}, \:
        (\sum_{i=1}^{d}\sum_{j=1}^{d} 
       c_{i}r_{j,i}^{k-1}w_{j}, \sum_{j=1}^{d}a_{k}^{j}w_{j})\in
       \mathbf{Z}\}=\\ 
     & =\{ w= \sum_{i=1}^{d}c_{i}w_{i}^{k-1}, \:
    \sum_{i=1}^{d}c_{i}(\sum_{j=1}^{d} a_{k}^{j}r_{j,i}^{k-1})\in
    \mathbf{Z}\}.
  \end{array}$$
We get:

\begin{lemm} \label{cal2}
  One has the equality of matrices:
   $$\mathfrak{m}( \mathcal{W}_{k-1}, \sigma_{k-1}, \prec_{k-1};
\mathcal{W}_{k})= \mathfrak{m}( \sum_{j=1}^{d}
a_{k}^{j}r_{j,1}^{k-1}, ..., \sum_{j=1}^{d}
a_{k}^{j}r_{j,d}^{k-1}).$$
\end{lemm}

The lemmas \ref{cal3} and \ref{cal2} allow to compute recursively the
matrices \linebreak  $\mathfrak{m}(\mathcal{W}, 
\sigma, \prec;\mathcal{W}_{k})$ for $k \in \{1,...,g\}$ from the
knowledge of the components of $a_{1},...,a_{k}$ in the basis
$(u_{1},...,u_{d})$. Indeed:
\begin{equation} \label{cal4}
 \mathfrak{m}(\mathcal{W},
 \sigma, \prec;\mathcal{W}_{k})= n_{s}(\mathfrak{m}(\mathcal{W}_{k-1},
 \sigma_{k-1}, \prec_{k-1};\mathcal{W}_{k})\mathfrak{m}(\mathcal{W},
 \sigma, \prec;\mathcal{W}_{k-1})).
\end{equation}

Once the matrix $\mathfrak{m}(\mathcal{W},
 \sigma, \prec;\mathcal{W}_{k})$ is known, lemma \ref{transf} shows
 that \linebreak $\mathfrak{m}(\mathcal{W}_{k},\sigma, \prec)$ is also known.

\medskip

In the special case in which $(\mathcal{S},0)$ is an irreducible
quasi-ordinary singularity of hypersurface having $(A_{1},...,A_{G})$
as characteristic exponents with respect to some projection, we put
$g=G, \:\mathcal{W}=W_{0}, \: \mathcal{M}=M_{0}$ and $a_{k}=A_{k}, \:
\forall \: k \in \{1,...,g\}$. By combining corollary \ref{normtor1}
and definition \ref{type}, we see that the normalization of
$(\mathcal{S},0)$ is a Hirzebruch-Jung singularity of type
$\mathfrak{m}(W_{G},\sigma_{0}, \prec_{0})$, which can be
computed by the previous method. Using the lemmas \ref{transf} (more
precisely the relation (\ref{normtrans})), \ref{cal3}, 
\ref{cal2} and relation (\ref{cal4}), we get the following compact
form of the algorithm: 

\begin{prop} \label{algfond}
  Let $f\in \mathbf{C}\{X_{1},...,X_{d}\}[Y]$ be an irreducible
  quasi-ordinary polynomial with characteristic exponents
  $A_{1},...,A_{G}$. We look at $A_{k}$ as a matrix $1 \times d$. If
  $R^{k}:= 
  \mathfrak{m}(W_{0},\sigma_{0}, \prec_{0};W_{k}), \:
  S^{k}:= \mathfrak{m}(W_{k-1},\sigma_{k-1},
  \prec_{k-1};W_{k}), \:
  T^{k}=(t_{i,j}^{k})_{i,j}:=(R^{k})^{-1}, \: d_{j}^{k}:=
  \mathrm{lcm}(\mathrm{den}(t_{1,j}^{k}),...,
  \mathrm{den}(t_{j,j}^{k}))), \: \forall \: k \in \{1,...,G\}, \:
  \forall \: j\in \{1,...,d\}$, and
  $R^{0}:=I_{d}$, then:
$$\begin{array}{l}
     S^{k}=
     \mathfrak{m}(A_{k}R^{k-1}),\\
         R^{k}= n_{s}(S^{k}R^{k-1})\\
      N_{k}=\mathrm{det}(S^{k}).
   \end{array}$$
The normalization of the germ defined by $f=0$ is a Hirzebruch-Jung 
singularity of type: 
$$\mathfrak{m}(W_{G},\sigma_{0}, \prec_{0})=
n_{a}((d_{j}^{G}t_{i,j}^{G})_{i,j}).$$ 
\end{prop}

We recall that the numbers $N_{k}$ were defined after the theorem
\ref{normtor}. 

If $(\mathcal{S},0)$ is an irreducible quasi-ordinary singularity of
dimension $d \geq 1$ and embedding dimension $d+1$, Lipman \cite{L 65}
showed that there is always a quasi-ordinary polynomial $f \in
\mathbf{C}\{X_{1},...,X_{d}\}[Y]$ defining $\mathcal{S}$ such that its
characteristic exponents $A_{1},...,A_{G}$ verify:
\begin{equation} \label{rel 6}
  \left\{ \begin{array}{l}
  (A_{1}^{1},...,A_{G}^{1})\geq_{lex}\cdots
  \geq_{lex}(A_{1}^{d},...,A_{G}^{d}) \\
     A_{1}^{2}\neq 0 \: \mbox{or} \: A_{1}^{1}>1
           \end{array} \right.
\end{equation}
Lipman \cite{L 88} and Gau \cite{G 88} showed that a sequence
$A_{1},...,A_{G}$ which verifies (\ref{rel 6}) - they called it then
\textit{normalized} - is an
embedded topological invariant of $(\mathcal{S},0)$. In particular, it
is an analytical invariant of $(\mathcal{S},0)$. In \cite{PP 03'} we gave
an algebraic proof of this analytical invariance. This shows that for
an irreducible quasi-ordinary germ of hypersurface, there is a way to
choose a well-defined matrix for the type of its normalization
between the $d!$ possibilities. Indeed, one simply starts the
application of the previous algorithm from normalized characteristic
exponents.

\section{A tridimensional example} \label{example}

Consider the following sequence of characteristic exponents:
$$A_{1}=(\frac{1}{4}, \frac{1}{6}, \frac{1}{6}), \: 
   A_{2}=(\frac{3}{8}, \frac{5}{12}, \frac{7}{12}).$$
As the relations (\ref{rel6}) are verified, it is a normalized
sequence (see the definition in the last paragraph of the previous section).

Let us apply the algorithm summarized in proposition \ref{algfond}:

$$\begin{array}{l}
   R^{1}= S^{1}=\left( \begin{array}{ccc}
                   4 & 2 & 2\\
                   0 & 3 & 2\\
                   0 & 0 & 1
                 \end{array} \right) \\
    N_{1}=\mathrm{det}(S^{1})=12 \\
  T^{1}= (R^{1})^{-1}= \dfrac{1}{2^{2}\cdot 3}
            \left( \begin{array}{ccc}
                   3 & -2 & -2\\
                   0 & 4 & -8\\
                   0 & 0 & 12
                 \end{array} \right)=\left( \begin{array}{ccc}
                    \dfrac{1}{4}& -\dfrac{1}{6} & -\dfrac{1}{6}\\
                   0 &  \dfrac{1}{3}&
                   -\dfrac{2}{3} \\
                   0 & 0 & 1\\
                \end{array} \right)\\
  d_{1}^{1}=4, \: d_{2}^{1}=6, \:
  d_{3}^{1}=6 \\
      (d_{j}^{1}t_{i,j}^{1})_{i,j}=
          \left( \begin{array}{ccc}
                   1 & -1 & -1\\
                   0 & 2 & -4\\
                   0 & 0 & 6
                 \end{array} \right)  \\ 
  \mathfrak{m}(W_{1},\sigma_{0},
  \prec_{0})=n_{a}((d_{j}^{1}t_{i,j}^{1})_{i,j})= 
         \left( \begin{array}{ccc}
                   1 & -1 & -1\\
                   0 & 2 & -4\\
                   0 & 0 & 6
                 \end{array} \right) \\
  A_{2}R^{1}=  \left( \begin{array}{ccc}
                   \dfrac{3}{2} & 2 & \dfrac{13}{6}
                 \end{array} \right)  \\
\end{array}$$
$$\begin{array}{l}
  S^{2}= \mathfrak{m}(\dfrac{3}{2}, 2, \dfrac{13}{6})= 
           \left( \begin{array}{ccc}
                   2 & 0 & 1\\
                   0 & 1 & 0\\
                   0 & 0 & 3
                 \end{array} \right)  \\
     N_{2}=\mathrm{det}(S^{2})=6\\
  S^{2}R^{1}= \left( \begin{array}{ccc}
                   8 & 4 & 5\\
                   0 & 3 & 2\\
                   0 & 0 & 3
                 \end{array} \right)  \\
  R^{2}= n_{s}( S^{2}R^{1})= \left( \begin{array}{ccc}
                   8 & 4 & 5\\
                   0 & 3 & 2\\
                   0 & 0 & 3
                 \end{array} \right)  \\
  T^{2}=(R^{2})^{-1}= \dfrac{1}{2^{3}\cdot 3^{2}}
            \left( \begin{array}{ccc}
                   9 & -12 & -7\\
                   0 & 24 & -16\\
                   0 & 0 & 24
                 \end{array} \right)=\left( \begin{array}{ccc}
                    \dfrac{1}{2^{3}}& -\dfrac{1}{2\cdot
                   3} & -\dfrac{7}{2^{3}\cdot 3^{2}}\\
                   0 &  \dfrac{1}{3}&
                   -\dfrac{2}{3^{2}} \\
                   0 & 0 & \dfrac{1}{ 3}
                 \end{array} \right) \\
 d_{1}^{2}=2^{3}, \: d_{2}^{2}=6, \:
  d_{3}^{2}=2^{3}\cdot 3^{2} \\
      (d_{j}^{2}t_{i,j}^{2})_{i,j}=
          \left( \begin{array}{ccc}
                   1 & -1 & -7\\
                   0 & 2 & -16\\
                   0 & 0 & 24
                 \end{array} \right)  \\
  \mathfrak{m}(W_{1},\sigma_{0},
  \prec_{0})=n_{a}((d_{j}^{2}t_{i,j}^{2})_{i,j})= 
         \left( \begin{array}{ccc}
                   1 & -1 & -7\\
                   0 & 2 & -16\\
                   0 & 0 & 24
                 \end{array} \right)
\end{array}$$

This shows that the normalization of a quasi-ordinary hypersurface 
singularity with one characteristic exponent $(\dfrac{1}{4}, 
\dfrac{1}{6}, \dfrac{1}{6})$ has a Hirzebruch-Jung singularity of type
$\left( \begin{array}{ccc} 
                   1 & -1 & -1\\
                   0 & 2 & -4\\
                   0 & 0 & 6
                 \end{array} \right) $ and the normalization of a
               quasi-ordinary hypersurface 
singularity with two characteristic exponents $(\dfrac{1}{4},
\dfrac{1}{6}, \dfrac{1}{6}), \: (\dfrac{3}{8}, \dfrac{5}{12},
\dfrac{7}{12})$ has a Hirzebruch-Jung singularity of type 
$\left( \begin{array}{ccc}
                   1 & -1 & -7\\
                   0 & 2 & -16\\
                   0 & 0 & 24
                 \end{array} \right)$.  

By the comments made at the end of the previous section, we have
obtained like this well-defined normal forms for the types of the
normalizations of the considered quasi-ordinary singularities.

\section{The classical (2-dimensional) Hirzebruch-Jung singularities}
  \label{class}

In this section we restrict to the case of surfaces and we compare our
definition of the type of a Hirzebruch-Jung singularity with the one
given in Barth, Peters, Van de Ven \cite{BPV 84}. For details on
Hirzebruch's work \cite{H 53}, one should consult Brieskorn 
\cite{B  00}. 

Let $p_{1},...,p_{r}$ be a sequence of integers, such
that $p_{i} \leq -2, \: \forall i\in \{1,...,r\}$ and $r\geq
1$. Let $X$ be a smooth complex analytical surface containing a reduced
divisor with normal crossings $C$ whose components $C_{1},...,C_{r}$
are projective lines with self-intersections $C_{i}^{2}=p_{i},
\: \forall i \in \{1,...,r\}$, and such that: 
$$C_{i}\cdot C_{j}=\left\{ \begin{array}{l}
                               1, \: \mbox{if } \mid i-j \mid=1 \\
                               0, \: \mbox{either}
                           \end{array} \right.$$                           
Such a couple $(X,C)$ always exists. The curve $C$ is called a
\textit{Hirzebruch-Jung string} with self-intersection numbers
$p_{i}$. 

Define also the coprime numbers $(n,q)\in \mathbf{N}^{2}$, $0<q<n$ by
the formula:
\begin{equation} \label{JHfrac}
 \frac{n}{q}=\mid p_{1}\mid -\cfrac{1}{\mid
  p_{2}\mid-\cfrac{1}{\cdots -\cfrac{1}{\mid
  p_{r}\mid}}}
\end{equation}

Then one has the following theorem ((5.1) in \cite{BPV 84}):

\begin{theor} \label{clJH} (Hirzebruch)
 If $C \subset X$ is a Hirzebruch-Jung string with self-intersection
 numbers $p_{i}$ satisfying the relation (\ref{JHfrac}), then
 the germ obtained 
 by contracting $C$ to a point is analytically isomorphic to the
 normalization of the germ at the origin of the surface with equation
 $Y^{n}=X_{1}X_{2}^{n-q}$. 
\end{theor}

This motivates the following definition given in \cite{BPV 84}:

\begin{defin} \label{cldef}
 A normal germ of surface is said to be \textbf{a Hirzebruch-Jung singularity
 of type $\mathcal{A}_{n,q}$} if it is analytically isomorphic with the
 normalization at the origin of the surface with equation
 $Y^{n}=X_{1}X_{2}^{n-q}$. 
\end{defin}

\noindent \textbf{Remark: } In \cite{AS 03}, Aroca and Snoussi showed
 more generally that any 
normal quasi-ordinary singularity (i.e. Hirzebruch-Jung singularity in
our terms) is the normalization of a complete intersection germ
defined by binomial equations. 
\medskip

Let us see the relation between the classical normal form of
definition \ref{cldef} and the one
we introduced in the definition \ref{type}. 

\begin{prop} \label{comptype}
 The Hirzebruch-Jung singularity
 of type $\mathcal{A}_{n,q}$ following definition \ref{cldef}  is of
 type $\left( \begin{array}{cc} 
               1 & -q \\
               0 & n
             \end{array} \right)$
 following definition \ref{type}. 
\end{prop}

\textbf{Proof :} 
The polynomial 
$Y^{n}-X_{1}X_{2}^{n-q}$ is  quasi-ordinary  with only one \linebreak
characteristic exponent $A_{1}=(\dfrac{1}{n}, 
1-\dfrac{q}{n})$. Applying lemma \ref{cal3}, with \linebreak 
$\mathfrak{m}(W_{0}, \sigma_{0}, \prec_{0};W_{1})=
   \left( \begin{array}{cc}
              r_{1,1} & r_{1,2}\\
              0 & r_{2,2}
          \end{array} \right)$ 
and 
$\mathfrak{B}(W_{0}, \sigma_{0},
\prec_{0};W_{1})=(w_{1}^{1}, w_{2}^{1})$, we get
$l_{1}= \mathrm{den}(\dfrac{1}{n})=n, \: l_{2}=
\mathrm{lcm}(\mathrm{den}(\dfrac{1}{n}),\mathrm{den}(1-\dfrac{q}{n}))=n,
\: r_{1,1}=l_{1}=n, \: r_{2,2}=\dfrac{l_{2}}{l_{1}}=1, \: 
 r_{1,2}=-(l_{1}A^{2}_{1}r_{2,2})(l_{1}A_{1}^{1})^{-1}$ in $\mathbf{Z}/
r_{1,1}\mathbf{Z}$, which implies $r_{1,2}=q$. So:
$$\left\{ \begin{array}{l}
              w_{1}^{1}= nw_{1}\\
              w_{2}^{1}= qw_{1}+w_{2}
           \end{array} \right.$$
which implies:
$$\left\{ \begin{array}{l}
              w_{1}= \dfrac{1}{n}w_{1}^{1}\\
              w_{2}= w_{2}^{1}-\dfrac{q}{n}w_{1}^{1}
           \end{array} \right.$$

Then, following the notations of proposition \ref{normcone}, if
$v_{1}, v_{2}$ are the minimal generators of $W_{1}$ situated on the
edges $\mathbf{R}_{+}w_{1}, \mathbf{R}_{+}w_{2}$ of $\sigma_{0}$, we
deduce:
\begin{equation} \label{cal5}
  \left\{ \begin{array}{l}
              v_{1}= w_{1}^{1}\\
              v_{2}= nw_{2}^{1}- qw_{1}^{1}
           \end{array} \right.
\end{equation}
This shows that:
$$\mathfrak{m}(W_{1}, \sigma_{0}, \prec_{0})=
      \left( \begin{array}{cc}
               1 & -q \\
               0 & n
             \end{array} \right).$$
The proposition is proved. \hfill $\Box$

\medskip

In dimension 2, there are only two possible choices of the ordering of
the edges of $\sigma$, and so only two matrices for the type of a
Hirzebruch-Jung singularity (see the comments following definition
\ref{type}). The following proposition which relates them was probably
first proved by Hirzebruch \cite{H 53}:

\begin{prop} \label{amb} (Hirzebruch)
  If two  Hirzebruch-Jung singularities of types $\mathcal{A}_{n,q}$
  and $\mathcal{A}_{n',q'}$ are isomorphic, then $n=n'$ and $(q=q' \:
  \mathrm{or} \: qq'\equiv 1 (\mathrm{mod} \: n))$.
\end{prop}

\textbf{Proof: } Hirzebruch proved this result by looking at the
minimal desingularisations of the singularities. Both have as
exceptional divisors 
Hirzebruch-Jung strings with the same 
sequences of self-intersection numbers, but possibly reversed. Then
one makes computations using formula (\ref{JHfrac}).

Here we give another proof, which uses the orbifold map instead of the
minimal desingularization one. As showed by theorem \ref{mthm}, the
couple $(\mathcal{W}, \sigma)$ is well-defined up to isomorphism by
the analytical structure of the singularity. If one chooses the reverse
order of $\prec_{0}$ in the previous computations, one gets:
$$\left\{ \begin{array}{l}
              v_{2}= e_{1}^{1}\\
              v_{1}= n'e_{2}^{1}- q'e_{1}^{1}
           \end{array} \right.,$$
where $(e_{1}^{1}, e_{2}^{1})$ is a basis of $W_{1}$. Combining these
relations with (\ref{cal5}), we get first $n=n'$, as both measure the
index $(W_{1}:\mathbf{Z}v_{1}+\mathbf{Z}v_{2})$. Then
$w_{1}^{1}=v_{1}=ne_{2}^{1}-
q'e_{1}^{1}=ne_{2}^{1}-q'v_{2}=ne_{2}^{1}-q'(nw_{2}^{1}-qw_{1}^{1}) \:
\Rightarrow \: (1-qq')w_{1}^{1}=n(e_{2}^{1}-q'w_{2}^{1})$. As
$w_{1}^{1}$ is a primitive element of $W_{1}$, this implies that
$1-qq' \equiv 0 (\mathrm{mod} \: n)$, which proves the proposition.

Another method would have been to apply the algorithm of normalization
as in the proof of proposition \ref{comptype}, but starting with
$A_{1}=(1-\dfrac{q}{n}, 
\dfrac{1}{n})$. \hfill $\Box$

\medskip

The computations we have done in order to prove proposition
\ref{comptype} are a particular case of the normalization algorithm
\ref{algfond} 
presented in the previous section. By using lemma \ref{cal3}, we can 
 give in a more explicit form
this algorithm, as we published it (but with slightly different
notations) in \cite{PP 01} and \cite{PP 03}:

\begin{prop} \label{2alg}
Let $f \in \mathbf{C}\{X_{1}, X_{2}\}[Y]$ be an irreducible
quasi-ordinary polynomial with characteristic exponents
$A_{1},...,A_{G}$. If $\mathfrak{m}(W_{0}, \sigma_{0},
\prec_{0};W_{k})=$ \linebreak
     $\left( \begin{array}{cc}
               r_{1,1}^{k} &  r_{1,2}^{k}\\
               0 & r_{2,2}^{k}
             \end{array} \right)$, 
$\mathfrak{m}(W_{k-1}, \sigma_{k-1},
\prec_{k-1};W_{k})=
     \left( \begin{array}{cc}
               s_{1,1}^{k} &  s_{1,2}^{k}\\
               0 & s_{2,2}^{k}
             \end{array} \right), \: \forall \: k \in \{1,...,G\}$,
and $\left( \begin{array}{cc}
               r_{1,1}^{0} &  r_{1,2}^{0}\\
               0 & r_{2,2}^{0}
             \end{array} \right)= 
    \left( \begin{array}{cc}
               1 &  0\\
               0 & 1
             \end{array} \right)$, then:
$$\begin{array}{l}
    s_{1,1}^{k} =\mathrm{den}(A_{k}^{1}r_{1,1}^{k})\\
    s_{2,2}^{k}
      =\dfrac{\mathrm{lcm}(\mathrm{den}(A_{k}^{1}r_{1,1}^{k}), 
         \mathrm{den}(A_{k}^{1}r_{1,2}^{k-1} + A_{k}^{2}r_{2,2}^{k-1}))
      }{\mathrm{den}(A_{k}^{1}r_{1,1}^{k})}\\ 
    s_{1,2}^{k}= \left\{ 
       \begin{array}{ll}
           -\mathrm{den}(A_{k}^{1}r_{1,2}^{k-1} +
            A_{k}^{2}r_{2,2}^{k-1})\cdot & \\
             \: \: \mathrm{lcm}(\mathrm{den}(A_{k}^{1}r_{1,1}^{k}),  
           \mathrm{den}(A_{k}^{1}r_{1,2}^{k-1} +
            A_{k}^{2}r_{2,2}^{k-1})) \cdot & \\
             \: \:       \mathrm{num}(A_{k}^{1}r_{1,1}^{k-1})^{-1} 
            \mbox{ in } \mathbf{Z}/
            \mathrm{den}(A_{k}^{1}r_{1,1}^{k-1})\mathbf{Z}, &  
            \mbox{ if }  \mathrm{den}(A_{k}^{1}r_{1,1}^{k-1}) \neq 1\\
              0, &  \mbox{ if }  \mathrm{den}(A_{k}^{1}r_{1,1}^{k-1})=1    
        \end{array}  \right.    \\
    r_{1,1}^{k} = s_{1,1}^{k}r_{1,1}^{k-1}\\
    r_{2,2}^{k} = s_{2,2}^{k}r_{2,2}^{k-1}\\
    r_{1,2}^{k} = s_{1,1}^{k}r_{1,2}^{k-1} +
    s_{2,2}^{k}r_{2,2}^{k-1} \mbox{ in } \mathbf{Z}/
            r_{1,1}^{k}\mathbf{Z}
  \end{array}$$
The normalization of the germ defined by $f=0$ is a Hirzebruch-Jung
singularity of type $\left( \begin{array}{cc}
               1 &
               -\dfrac{r_{1,2}^{G}}{\mathrm{gcd}(r_{1,2}^{G},r_{1,2}^{G})} \\ 
               0 & \dfrac{r_{1,1}^{G}}{\mathrm{gcd}(r_{1,2}^{G},r_{1,2}^{G})}
             \end{array} \right)$. 
\end{prop}

\section{Questions} \label{quest}

If $(\mathcal{S},0)$ is a reduced germ of complex analytical space, we
denote by $K(\mathcal{S})$ its \textit{abstract boundary}. It is defined as the
intersection of a representative of $(\mathcal{S},0)$ with a
sufficiently small euclidean sphere centered at $0$ in an arbitrary
system of local coordinates at $0$. It is independent of these choices
(Durfee's proof in \cite{D 83} for algebraic varieties extends to
analytical ones).

Hirzebruch \cite{H 53} noticed that the abstract boundary of a
bidimensional Hirze\-bruch-Jung singularity $(\mathcal{Z},0)$ of type
$\mathcal{A}_{n,q}$ is a lens space $L(n,q)$. As it was known since
Reidemeister \cite{R 35} that $L(n,q)$ is homeomorphic to
$L(n',q')$ if and only if $n=n'$ and $(q=q' \: \mathrm{or} \:
qq'\equiv 1 (\mathrm{mod} \: n))$, this showed by proposition
\ref{amb} that in this case the
homeomorphism type of $K(\mathcal{Z})$ determines the analytical type of
$(\mathcal{Z},0)$. More generally, we ask:

\medskip

\textbf{Question 1} \textit{Let  $(\mathcal{Z},0)$ be a
  Hirzebruch-Jung singularity of dimension $\geq 3$. Denote by
  $K(\mathcal{Z})$ its abstract 
  boundary. Is it true that the homeomorphism type of $K(\mathcal{Z})$
  determines the analytical type of $(\mathcal{Z},0)$?}

\medskip

If the answer to the previous question is negative, it would be
interesting to know what supplementary structure one should add to the
boundary  $K(\mathcal{Z})$ in order to make it affirmative (e.g. should one
consider it rather as an orbifold, or add some stratified smooth structure?)

In the case when the canonical representation $\rho(\mathcal{Z})$
associated to the singularity (see (\ref{repr})) 
is a cyclic fixed-point free action outside the origin, the answer to
the question is affirmative. Indeed, in this case the boundary is a 
\textit{generalized lens space} and the corresponding result was
obtained by Franz 
(see Dieudonn{\'e} \cite{D 89}, page 246). If $d \geq 3$, the action
$\rho(\mathcal{Z})$ 
may be non-cyclic, and even if it is cyclic, it may have fixed
points. One can decide if $\Gamma(\mathcal{Z})$ is cyclic by computing the
invariant factors of a matrix of presentation of $\widetilde{\mathcal{W}}$
with respect to $\mathcal{W}$, for example $\mathfrak{m}(\mathcal{W},
\sigma, \prec)$ for an arbitrary ordering $\prec$.
\medskip

Consider now more general pairs $(\mathcal{W}, \sigma)$ than the
simplicial ones:

\medskip

\textbf{Question 2} \textit{Let $(\mathcal{W}, \sigma)$ be a maximal
  pair, where $\sigma$ is not a simplicial cone. If $0$ denotes the
  0-dimensional orbit of the affine toric variety
  $\mathcal{Z}(\mathcal{W}, \sigma)$, is it true that the analytic
  type of the germ $(\mathcal{Z}(\mathcal{W}, \sigma),0)$ determines
  the pair $(\mathcal{W}, \sigma)$ up to isomorphism?}
\medskip

In this case, the  germ $(\mathcal{Z},0)$ is not a  finite quotient
singularity. So, in order to attack this question, it seems
that one cannot avoid anymore the use of some desingularization
morphism. A first step towards the solution could come from an 
affirmative answer to the following question:
\medskip

\noindent
\textbf{Question 3} \textit{Could one prove theorem \ref{isom} using
resolutions of the singularities of the germ $(\mathcal{Z},0)$ instead
of the canonical 
representation $\rho(\mathcal{Z})$ of its local fundamental group?}

\medskip

By analogy with question 1, we ask also:
\medskip

\noindent
\textbf{Question 4} \textit{Let 
$\mathcal{Z}$ be an affine (not necessarily simplicial) toric
variety. Is it true that the homeomorphism type (possibly enriched
with supplementary structure) of
$K(\mathcal{Z})$ determines the analytical type of $(\mathcal{Z},0)$?}

\medskip

\noindent Patrick Popescu-Pampu \\
Univ. Paris 7 Denis Diderot \hspace{3mm} \\Instit. de
Maths. - UMR CNRS 7586 \\ Equipe "G{\'e}om{\'e}trie et dynamique" \hspace{3mm}\\
Case 7012\\ 
2, place Jussieu  \hspace{3mm}\\ 75251-Paris Cedex 05 \hspace{3mm}\\ FRANCE
\medskip

\end{document}